\documentclass[11pt]{amsart}
\usepackage[english]{babel}
\usepackage{physics}
\usepackage{braket}
\usepackage{dsfont}
\usepackage[table,xcdraw]{xcolor}
\usepackage{multirow}
\usepackage{graphicx}
\usepackage{tikz}
\usepackage{amssymb}
\usepackage{mathtools}
\usetikzlibrary{matrix}

\usepackage[margin=1in]{geometry}
\usepackage[pagebackref, colorlinks = true, linkcolor = blue, urlcolor  = blue, citecolor = red]{hyperref}
\usepackage[utf8x]{inputenc}

\renewcommand{\epsilon}{\varepsilon}

\newtheorem*{thm*}{Theorem}

\newtheorem*{specthm2*}{Theorem 7}
\newtheorem*{specthm*}{Theorem 7 revised}

\mathtoolsset{showonlyrefs}

\title{Comment on ``Generators of matrix algebras in dimension 2 and 3"}

\author{{\'A}ngela Capel}
\email{angela.capel@ma.tum.de}
\address{Department of Mathematics, Technische Universit\"at M\"unchen, 85748 Garching, Germany and Munich Center for Quantum Science and Technology (MCQST), M\"unchen, Germany}

\author{Yifan Jia}
\email{yifan.jia@tum.de}
\address{Department of Mathematics, Technische Universit\"at M\"unchen, 85748 Garching, Germany and Munich Center for Quantum Science and Technology (MCQST), M\"unchen, Germany}
\date{ \today. \\
\phantom{asl}\textit{Keywords:} Generator; Matrix; Algebra. \\
\phantom{asl}\textit{2010 Mathematics Subject Classification:}  Primary 15A30; Secondary 47L05.}

\begin{document}

\maketitle

\begin{abstract}
    Theorem 7 in \cite{AslaksenSletsjoe_WrongPaper_2009} states sufficient conditions to determine whether a pair generates the algebra of $3 \times 3$ matrices over an algebraically closed field of characteristic zero. In that case, an explicit basis for the full algebra is provided, which is composed of words of small length on such pair. However, we show that this theorem is wrong, since it is based on the validity of an identity which is not true in general.
    
\end{abstract}

\section{discussion}

Let $M_n (\mathbb{K})$ denote the set of all $n \times n $ matrices over a field $\mathbb{K}$. Let $S$ be a subset of  $M_n (\mathbb{K})$ and denote by $S^m$ the set of all products of the form $A_1 \cdots A_m$, with $A_i \in S\cup\{I_n\}$ for all $i=1, \ldots , m$, where $I_n$ is the $n \times n$ identity matrix. We say that a generating set $S$ has length $k \in \mathbb{N}$ if 
\begin{equation*}
    \text{span}\left\lbrace S^{k}  \right\rbrace = M_n (\mathbb{K}) \, \, , \phantom{sdad}\text{and } \, \, \; \; \; \text{ span}\left\lbrace S^{k-1}  \right\rbrace \subsetneq \text{span}\left\lbrace S^{k}  \right\rbrace \, .
\end{equation*}
The problem of finding bounds on the length of generating sets, and in particular generating pairs, has been thoroughly studied in the past decades. For arbitrary order $n$, the best known bound on the length of any generating set is $O(n\log n)$ \cite{shitov2019improved}, although it was conjectured years before that the optimal bound might be $2n-2$ \cite{Paz_Conjecture_1984}. Indeed, this is the case at least for $n \leq 6$ \cite{LongstaffNiemeyerPanaia_LengthDimension5_2006, LambrouLongstaff_LengthDimension6_2009}, for which the bound is shown to be sharp. Moreover, this bound also holds for arbitrary $n$ under different conditions on one of the generators \cite{LogstaffRosenthal_LengthsIrreduciblePairs_2011, GutermanLaffeyMarkovaSmigoc_PazConjecture_2018}, even though it is always possible to find $n$ matrices in  $M_n (\mathbb{K})$ such that the words of length $2$ in those matrices span the whole algebra \cite{Rosenthal_WordBasis_2012}. Finally, when the problem is reduced to the study of generic matrices, the bound can be arguably improved to $O(\log n)$ \cite{KlepSpenko_SweepingWords_2016}.

In \cite{AslaksenSletsjoe_WrongPaper_2009}, the problem of providing conditions under which a set of $2 \times 2$ or $3 \times 3$ matrices over an algebraically closed field of characteristic zero generate the full matrix algebra is addressed. However, Theorem 7 in that paper is wrong, since Equation (3) in that reference does not hold in general. In this note, we give some comments on that result and present numerical evidence to show the falseness of the aforementioned expression  for arbitrary matrices $A$ and $B$. Unless stated otherwise, we follow  \cite{AslaksenSletsjoe_WrongPaper_2009} for terminology and notations. 

Before recalling Theorem 7 of \cite{AslaksenSletsjoe_WrongPaper_2009}, we denote by $ M_3$ the algebra of $3 \times 3$ matrices over a certain algebraically closed field of characteristic zero. Moreover, we write $[A,B]$ for the commutator of two matrices $A$ and $B$ and we define 
\begin{equation*}
    H(M):= \frac{ \tr[M]^2 - \tr[M^2]}{2} \, ,
\end{equation*}
where $\tr[M]$ denotes the trace of a matrix $M$. We can now state the aforementioned result, namely Theorem 7 of \cite{AslaksenSletsjoe_WrongPaper_2009}.

\newpage

\begin{specthm2*}\label{thm:wrong}
   Let $A,B\in M_3$. Then 
   \begin{align}\label{eq:false-expression}
       \det(I,A,A^2,B,B^2,AB,BA,[A,[A,B]],[B,[B,A]])=9\det [A,B]H([A,B]),
   \end{align}
   so if $\det [A,B]\neq0$ and $H([A,B])\neq 0$, then 
   \begin{align}
       \{I,A,A^2,B,B^2,AB,BA,[A,[A,B]],[B,[B,A]]\}
   \end{align}
   form a basis for $M_3$.
\end{specthm2*}

Note that this result intends to provide a sufficient condition for a pair of matrices $\left\lbrace A, B \right\rbrace$ to generate the full algebra $M_3$, as well as construct a basis for the algebra from words in such matrices. In particular, if the result was correct, it would yield the fact that for pairs $\left\lbrace A, B \right\rbrace$ such that $\det [A,B]\neq0$ and $H([A,B])\neq 0$, the necessary length of words to generate the full $M_3$ is 3, improving thus for this subclass of pairs the bound of $2n-2$ obtained in \cite{LongstaffNiemeyerPanaia_LengthDimension5_2006}. Moreover, such matrices are in particular generic, so they belong to the case studied in \cite{KlepSpenko_SweepingWords_2016}, in which the tighter bound for the length of words in dimension $3$ is $2 \lceil \log_2 (3) \rceil =4$. Thus, Theorem 7 would also provide an improvement to the latter result in dimension $3$.

\begin{figure}[h!]\label{fig:numericaltest}
    \centering
    \includegraphics[scale=0.7]{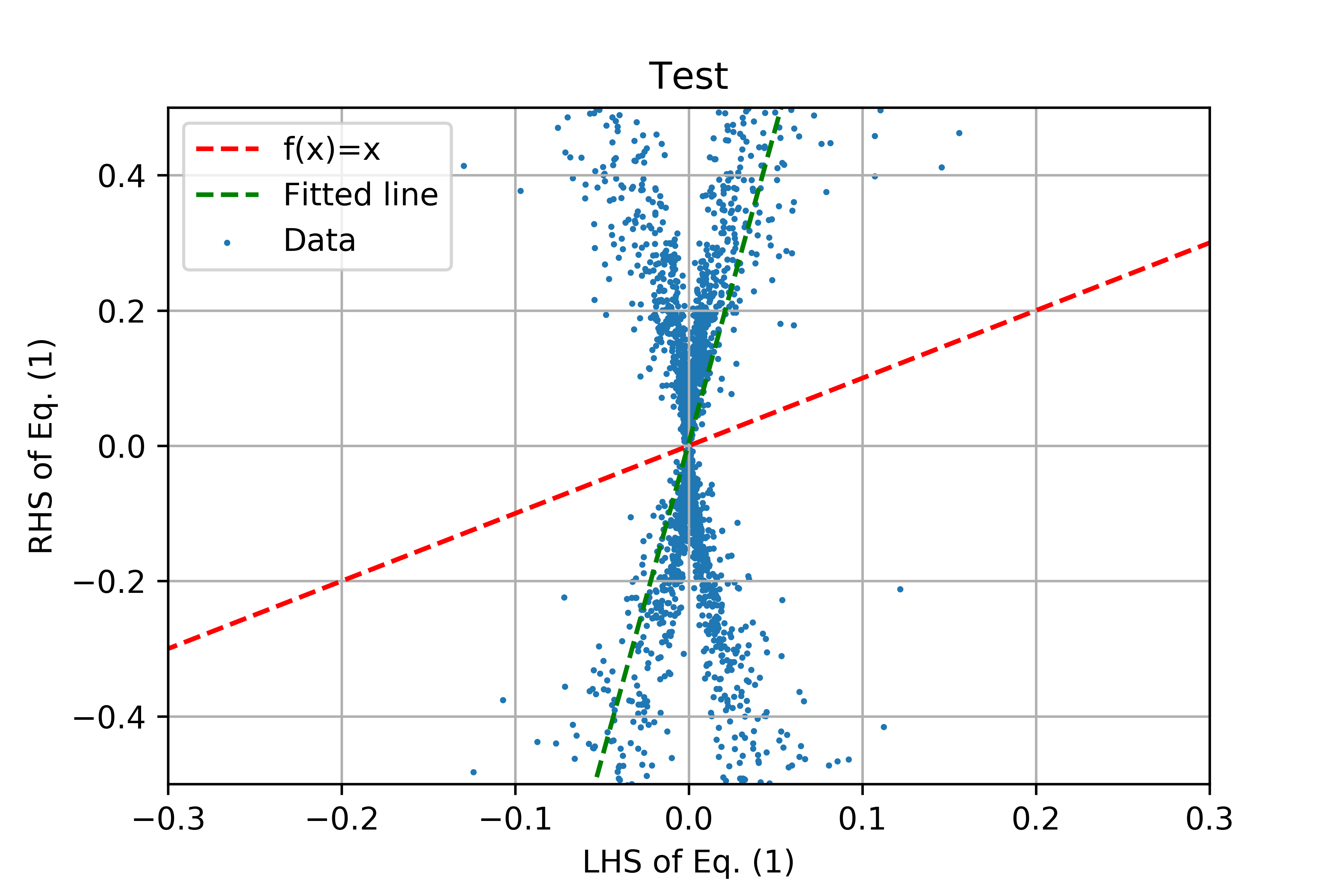}
    \caption{Test for 5000 randomly chosen $3\times3$ real matrix pairs. There exists no line such that every point fits to it, leading to the result that the factorization cannot be exact. However, there exists a strong correlation between the absolute values of both sides of Equation \eqref{eq:false-expression}. }
    \label{commenttest}
\end{figure}

To show that Equation \eqref{eq:false-expression} is false, and thus Theorem 7 does not hold in general, it is enough to construct a counterexample. For any $\delta >0 , \,  \varepsilon>0$, consider the following pair of matrices:
\begin{align}
    A=\begin{pmatrix}
     1&0&1\\0&-1&\delta\\1&0&1
      \end{pmatrix},\quad
      B=\begin{pmatrix}
     -1&0&1\\0&1&0\\1&\varepsilon&-1
      \end{pmatrix}.
\end{align}
Then, a short computation yields the following:
\begin{align}
   &\det(I,A,A^2,B,B^2,AB,BA,[A,[A,B]],[B,[B,A]])=-27 \,  \delta^5 \varepsilon^5+9 \, \delta^6 \varepsilon^6 \, ,\\
   &9\det [A,B]H([A,B])=27 \, \delta^3 \varepsilon^3- 9 \, \delta^4 \varepsilon^4.
\end{align}
Clearly, we can conclude that, for $A$ and $B$ as above, the left- and right-hand sides of Equation \eqref{eq:false-expression} are different, although we appreciate that both vanish if, and only if, either $\varepsilon$ or $\delta$ is zero. Thus, despite being different, both sides of the equation present a strong correlation.  This is noticeable in Figure \ref{commenttest}.

In Figure \ref{fig:numericaltest}, we generate 5000 random matrix pairs and compare their value when we insert the pair into the LHS and the RHS of Equation \eqref{eq:false-expression}. The points do not lie on some line $f(x)=ax+b$, especially not on $f(x)=x$. Therefore, we exclude the possibility that the incorrectness is caused by the improper calculation of coefficient $9$ in the RHS of Equation \eqref{eq:false-expression}. Moreover, if we restrict the data to the first and the third quadrants close to zero, the best fitted line has  slope $a\approx9$ and  intercept $b\approx 0$. For example, for the data in Figure \ref{fig:numericaltest}, the program returns $f(x)\approx9.3138x+0.0044$ for the fitted line, illustrated by the green line. If we consider the data that is further from the origin, the absolute value of the slope decreases and the points are not so concentrated as the image near 0, which is similar to the shape $``x"$.    

To conclude, we have shown that Theorem 7 in \cite{AslaksenSletsjoe_WrongPaper_2009} is false. For a correct upper bound on the necessary length of words on a pair to generate $M_3$, we refer the reader to  \cite{LongstaffNiemeyerPanaia_LengthDimension5_2006} for the general case and to \cite{KlepSpenko_SweepingWords_2016} for the generic case. 

\vspace{0.3cm}
\noindent
\hrulefill

\vspace{0.3cm}

After the completion of this note, we realized that for Equation \eqref{eq:false-expression} to be true, there were a factor $\det [A,B]$ and a minus sign missing in the right hand side. Therefore, we present the correct form of Theorem 7 below.

\begin{specthm*}\label{thm:true}
   Let $A,B\in M_3$. Then 
   \begin{align}\label{eq:true-expression}
       \det(I,A,A^2,B,B^2,AB,BA,[A,[A,B]],[B,[B,A]])=-9(\det [A,B])^2 H([A,B]),
   \end{align}
   so if $\det [A,B]\neq0$ and $H([A,B])\neq 0$, then 
   \begin{align}
       \{I,A,A^2,B,B^2,AB,BA,[A,[A,B]],[B,[B,A]]\}
   \end{align}
   form a basis for $M_3$.
\end{specthm*}

This yields the fact that whenever $\det [A,B]\neq0$ and $H([A,B])\neq 0$, the set composed of words on matrices $A, B$, which was presented in the original paper \cite{AslaksenSletsjoe_WrongPaper_2009}, actually spans the whole matrix algebra $M_3$. In particular, the necessary length of words to generate the full $M_3$ is 3. As mentioned above, this indeed improves the bound  $2n-2$ obtained in \cite{LongstaffNiemeyerPanaia_LengthDimension5_2006}, as well as the bound $2 \lceil \log_2 (n) \rceil$ provided in \cite{KlepSpenko_SweepingWords_2016} for generic matrices in dimension $n=3$.

\vspace{0.3cm}

\noindent {\it Acknowledgments.} 
We thank Michael Wolf for his comments and suggestions. This work has been partially supported by the Deutsche Forschungsgemeinschaft (DFG, German Research Foundation) under Germany's Excellence Strategy EXC-2111 390814868. YJ acknowledges support from the TopMath Graduate Center of the TUM Graduate School and the TopMath Program of the Elite Network of Bavaria.

\vspace{0.3cm}

\bibliographystyle{alpha}
\bibliography{biblio}

\vspace{0.3cm}

\end{document}